\documentclass{article}

\usepackage{amsmath,amsthm,amssymb,amscd,mathabx,authblk}
\usepackage{amsfonts}
\usepackage{authblk}
\usepackage{rotating,enumerate,enumitem,mathtools}
\usepackage{euscript}
\usepackage{pst-node}
\usepackage{epsfig,verbatim}
\usepackage[symbol]{footmisc}

\def\be{\begin{equation}}
\def\ee{\end{equation}}

\def\C{{\mathbb C}} 
\def\f{\EuScript}
 
\def\P{{\mathbb P}}
\def\Z{{\mathbb Z}}

\def\e{\eqref}
\def\phi{{\varphi}}
 
\def\tt{\widetilde}
\def\deg{{\rm deg\,}}

\def\bp{\begin{proposition}}
\def\ep{\end{proposition}}

\def\bt{\begin{theorem}}
\def\et{\end{theorem}}
\def\br{\begin{remark}}
\def\er{\end{remark}}
\def\be{\begin{equation}}
\def\bee{\begin{equation*}}
\def\l{\label}

\def\e{\eqref}

\def\ee{\end{equation}}
\def\eee{\end{equation*}}
\def\bl{\begin{lemma}}
\def\el{\end{lemma}}
\def\bp{\begin{problem}}
\def\ep{\end{problem}}
\def\bc{\begin{corollary}}
\def\ec{\end{corollary}}
\def\pr{\noindent{\it Proof. }}

\def\bd{\begin{definition}}
\def\ed{\end{definition}}
\def\t{\widetilde}

\mathtoolsset{showonlyrefs}

\newtheorem{theorem}{Theorem}[section]
\newtheorem{lemma}{Lemma}[section]
\newtheorem{definition}{Definition}[section]
\newtheorem{corollary}{Corollary}[section]
\newtheorem{proposition}{Proposition}[section]
\newtheorem{remark}{Remark}[section]
\newtheorem{problem}{Problem}[section]

\renewcommand{\thefootnote}{\fnsymbol{footnote}}

\begin{document}

\title{On mutually semiconjugate rational functions}
\author{F. Pakovich
\thanks{Department of Mathematics, Ben Gurion University of the Negev, Israel. e-mail:
pakovich@math.bgu.ac.il}
}

{
\makeatletter
\addtocounter{footnote}{1} 
\renewcommand\thefootnote{\@fnsymbol\c@footnote}%
\makeatother
\maketitle
}

\maketitle

\begin{abstract} We characterize pairs of rational functions $A$, $B$ such that 
$A$ is semiconjugate to $B$, and $B$ is semiconjugate to $A$.
\end{abstract} 

\begin{section}{Introduction}
Let $A$ and $B$ be rational functions of degree at least two on the Riemann sphere. 
The function $B$ is said to be {\it semiconjugate} to the function $A$ if there exists a non-constant rational function $X$
such that the diagram
\be \l{i1}
\begin{CD}
\C\P^1 @>B>> \C\P^1\\
@VV X V @VV X V\\ 
\C\P^1 @>A >> \C\P^1\ 
\end{CD}
\ee
commutes. 
 If $X$ is invertible, the functions $A$ and $B$ are called {\it conjugate}.  
In terms of dynamical systems, the conjugacy condition means that the dynamical systems $A^{\circ k},$ $k\geq 1,$ and $B^{\circ k},$ $k\geq 1,$  on $\C\P^1$ are equivalent, while the more general condition 
\eqref{i1} means that the first of these systems  is a   factor of the second. In particular, \eqref{i1} implies that $X$ sends attracting, repelling, and indifferent periodic 
points of $B$ to periodic points of $A$ of the same character. 
Note that the semiconjugacy relation is not symmetric. However, it is clear that if  $B$ is semiconjugate to $A$, and  $C$ is semiconjugate to $B,$ then $C$  is semiconjugate to $A.$ 
 Therefore, the
semiconjugacy relation is a preorder on the set of rational functions.

Although semiconjugate rational functions appear naturally in complex and arithmetic  dynamics (see e.g. the  papers \cite{be},  \cite{e},  \cite{ms}, \cite{pj}),  the problem of describing such functions started to be systematically studied only recently 
in the series of papers \cite{semi},  \cite{arn}, \cite{rec},   \cite{gl}, \cite{pnew}.
In this paper we address the following related question: under what conditions rational functions 
$A$ and $B$ are {\it mutually} semiconjugate, that is  $A$ is semiconjugate to $B$, and $B$ is semiconjugate to $A$ ? Such functions are of interest since
they exhibit very similar although not identical dynamics. In fact,  the mutual semiconjugacy relation can be considered  as a weaker form of the classical conjugacy relation. 

Examples of mutually semiconjugate rational functions can be obtained by the 
following construction. Let $A$ be a rational function.   
For any decomposition $A=U\circ V$ of $A$ into a composition 
of rational functions, we say that the 
rational function $\t A=V\circ U$ is  an {\it elementary transformation} of $A$. We say that rational functions $A$ and $B$ are {\it equivalent} and write $A\sim B$ if there exists 
a chain of elementary transformations between $A$ and $B$. 
Since obviously 
$$\t A\circ V=V\circ A, \ \ \  \ \ \ A\circ U=U\circ \t A,$$  
elementary transformations are mutually semiconjugate,
implying inductively that functions  $A$ and $B$  are mutually semiconjugate whenever $A\sim B$. Moreover, the corresponding semiconjugacy map $X$  preserves not only the  character of periodic points but also their exact periods and multipliers (see \cite{rec}).

Roughly speaking, the main result of this paper states that rational functions $A$ and $B$ are  mutually semiconjugate {\it only}
if $A\sim B$, unless $A$ and $B$ belong to the class of {\it Latt\`es maps},  which is  known to be a source of  exceptional 
examples in complex dynamics.  A typical example $A_{n,L}$ of such a map is   obtained 
from the ``multiplication theorem'' for the Weierstrass function: 
$$\wp_L(nz)=A_{n,L}\circ \wp_L(z),$$ 
where $\wp_L$ is the Weierstrass function with period lattice $L$, and $n\geq 2$ is an integer.
More precisely, we show that if mutually semiconjugate
rational functions $A$ and $B$ are not equivalent, then they are Latt\`es maps  with orbifold signature $(2,2,2,2)$.

\bt \l{main} Let $A$ and $B$ be mutually semiconjugate rational functions of degree at least two.
Then 
either $A\sim B$, or there exist orbifolds $\f O_1$ and $\f O_2$ with signature $(2,2,2,2)$ on the Riemann sphere such that 
$A:\f O_1\rightarrow \f O_1$ and $B:\f O_2\rightarrow \f O_2$ are covering maps between orbifolds. 
\et

Theorem \ref{main} implies that, apart from the very special class of Latt\`es maps, the equivalence relation induced by the mutual semiconjugacy coincides with the equivalence $\sim$ defined above. In particular, for a rational function $A$ that is not a Latt\`es map  there 
exist at most finitely many conjugacy classes of rational functions mutually semiconjugate to $A$, since a similar statement is true  for equivalence classes of $\sim$ (see \cite{rec}). 

The paper is organized as follows. In the second section, we recall some definitions and results  
concerning Riemann surface orbifolds, Latt\`es maps, and  commuting rational functions. We also prove a result concerning 
mutually semiconjugate Latt\`es maps with signatures distinct from $(2,2,2,2)$. In the third section, we review results about the equivalence $\sim$ and semiconjugate rational functions, and prove 
Theorem \ref{main}. 
Finally, in the fourth section we consider mutually semiconjugate Latt\`es maps with orbifold signature $(2,2,2,2)$, and 
construct examples of such maps that are not equivalent.  

\end{section}

\begin{section}{Orbifolds and commuting functions}

The problem of describing mutually semiconjugate rational functions is closely related to the problem of describing {\it commuting} rational functions. Indeed, if  $A$ and $B$ are mutually semiconjugate rational functions, then  there exist rational functions $X$ and $Y$ such that the diagram
\be \l{en}
\begin{CD}
\C\P^1 @>A>> \C\P^1\\
@VV Y V @VV Y V\\ 
\C\P^1 @>B >> \C\P^1\\ 
@VV X V @VV X V\\ 
\C\P^1 @>A >> \C\P^1\ 
\end{CD}
\ee
commutes, implying that the rational function $X\circ Y$ commutes with $A$. Similarly, the rational function $Y\circ X$ commutes with $B$. Commuting rational functions 
were investigated already by Julia \cite{j}, Fatou \cite{f}, and Ritt  \cite{r}. 
The most complete result, obtained by Ritt, states roughly speaking that commuting rational functions  
{\it having no iterate in common}
reduce either to powers, or to Chebyshev polynomials,  
or to Latt\`es maps.
 A  proof of the Ritt theorem based on modern dynamical methods was given by 
Eremenko \cite{e2}. Rational functions which {\it do} have  a common iterate were studied in \cite{rev}. 

In this paper, we will use the Ritt theorem in its modern formulation, given in \cite{e2}. This formulation uses the notion of {\it orbifold}.
Recall that a Riemann surface orbifold is a pair $\f O=(R,\nu)$ consisting of a Riemann surface $R$ together with a ramification function $\nu:R\rightarrow \mathbb N\cup \{\infty\}$ which takes the value $\nu(z)=1$ except at isolated points. 
For an orbifold $\f O=(R,\nu)$, 
 the {\it  Euler characteristic} of $\f O$ is the number
\be \l{euler} \chi(\f O)=\chi(R)+\sum_{z\in R}\left(\frac{1}{\nu(z)}-1\right),\ee
the set of {\it singular points} of $\f O$ is the set 
$$c(\f O)=\{z_1,z_2, \dots, z_s, \dots \}=\{z\in R \mid \nu(z)>1\},$$ and  the {\it signature} of $\f O$ is the set 
$$\nu(\f O)=\{\nu(z_1),\nu(z_2), \dots , \nu(z_s), \dots \}.$$ 
This definition of orbifold  (see e. g. \cite{dh}, \cite{e2}) differs slightly from  the definition given in \cite{mil}, say, where it is assumed that $\nu$ takes only {\it finite} values. To pass from the first definition to the second it is necessary to change the surface $R$ in the definition  of $\f O$ removing all points $z$ where $\nu(z)=\infty$. 
The same remark concerns other related definitions given below. 
Note that since removing a point from a surface $R$ reduces the Euler characteristic $\chi(R)$ by one, this passage does not change the Euler characteristic $\chi(\f O)$ defined by \eqref{euler}.

If 
 $R_1$, $R_2$ are Riemann surfaces provided with ramification functions $\nu_1,$ $\nu_2$, then a holomorphic branched covering map
$$f:\, R_1\setminus\{z:\nu_1(z)=\infty\} \rightarrow R_2\setminus\{z:\nu_1(z)=\infty\}$$ 
is called  {\it a covering map} $f:\,  \f O_1\rightarrow \f O_2$
between orbifolds
$\f O_1=(R_1,\nu_1)$ and $\f O_2=(R_2,\nu_2)$
if for any $z\in R_1$ the equality 
\be \l{us} \nu_{2}(f(z))=\nu_{1}(z)\deg_zf\ee holds.
It follows from the chain rule  that if 
$f:\, \f O_1\rightarrow \f O_2 $ and $g:\, \f O_2\rightarrow \f O_3$ are  covering maps between orbifolds, then $g\circ f:\, \f O_1\rightarrow \f O_3$ is also a co\-vering map. 
If $f:\,  \f O_1\rightarrow \f O_2$ is a covering map of finite degree between orbifolds with compact $R_1$ and $R_2$, then  the Riemann-Hurwitz 
formula implies that 
\be \l{rhor} \chi(\f O_1)=\chi(\f O_2)\deg f . \ee

{\it A universal covering} of an orbifold ${\f O}$
is a covering map between orbifolds  $\theta_{\f O}:\,
\tt {\f O}\rightarrow \f O$ such that $\tt R$ is simply connected and $\tt \nu(z)\equiv 1.$ 
If $\theta_{\f O}$ is such a map, then 
there exists a group $\Gamma_{\f O}$ of conformal automorphisms of $\tt R$ such that the equality 
$\theta_{\f O}(z_1)=\theta_{\f O}(z_2)$ holds for $z_1,z_2\in \tt R$ if and only if $z_1=\sigma(z_2)$ for some $\sigma\in \Gamma_{\f O}.$ A universal covering exists and 
is unique up to a conformal isomorphism of $\tt R,$
unless $\f O$ is the Riemann sphere with one ramified point or with two ramified points $z_1,$ $z_2$ such that $\nu(z_1)\neq \nu(z_2).$   
 Furthermore, 
$\tt R=\mathbb D$ if and only if $\chi(\f O)<0,$ $\tt R=\C$ if and only if $\chi(\f O)=0,$ and $\tt R=\C\P^1$ if and only if $\chi(\f O)>0$. Any covering map  $f:\,  \f O_1\rightarrow \f O_2$ between orbifolds lifts to an
isomorphism $\phi:\, 
\tt{R}_1\rightarrow \tt{R}_2$ which makes the diagram  
\be \l{esz}
\begin{CD}
\tt {R}_1 @>\phi >> \tt {R}_2\\
@VV\theta_{\f O_1}V @VV\theta_{\f O_2}V\\ 
\f O_1 @>f >> \f O_2\ 
\end{CD}
\ee
commutative, and maps points that are in the same orbit of $\Gamma_{\f O_1}$ to points that  are in the same orbit of $\Gamma_{\f O_2}$. The isomorphism $\phi$ is defined up to
a transformation $\phi\rightarrow g\circ \phi,$ where $g\in \Gamma_{\f O_2}$.
 In the other direction, for
any  isomorphism $\phi$ which maps any orbit of $\Gamma_{\f O_1}$ to an orbit of  $\Gamma_{\f O_2}$ 
there exists a uniquely defined covering map between orbifolds $f:\f O_1 \rightarrow \f O_2$ such that diagram \eqref{esz} commutes (see \cite{mil}, Appendix E,  and \cite{semi}, Section 3).  

Commuting rational functions having no iterate in common can be described in terms of  orbifolds 
$\f O=(\C\P^1,\nu)$ with $\chi(\f O)=0$. 
The signature of such an orbifold has one of the following forms
 \be \l{list} (\infty,\infty), \ \ \ (2,2,\infty), \ \ \  (2,2,2,2), \ \ \ (3,3,3), \ \ \ (2,4,4), \ \ \ (2,3,6).\ee
Correspondingly, the group $\Gamma_{\f O}$ is conjugate in $Aut(\C)$ to 
\begin{equation}\l{lis}
\begin{split} 
     &z\rightarrow    z + im,\ m\in \Z ; \\
  &z\rightarrow    \pm z + m,\ m \in \Z ; \\
     &z\rightarrow    \pm z + m + n\tau, \ m, n\in \Z ; \\
				&z\rightarrow    \omega^{2k}z+ m + n\omega,\quad m,n\in \Z, \ \ \  0\leq k \leq 2 ; \\
    &z\rightarrow    i^kz+ m +ni,\quad m, n\in   \Z, \ \ \ 0\leq k \leq 3 ; \\
      &z\rightarrow    \omega^{k}z+ m +n\omega, \quad    m,n\in \Z, \ \ \ 0\leq k \leq 5 ,
\end{split}
\end{equation}
where  $\tau$ is a complex number with  $\Im(\tau)>0$, and $\omega= e^{\pi i/3}$.
Finally, the universal covering of $\f O$ with $\Gamma_{\f O}$ from the list \eqref{lis}, up to the 
transformation 
$ \theta_{\f O}\rightarrow \mu \circ \theta_{\f O},$ where $\mu \in Aut(\C\P^1),$
 is
$$ exp(2\pi z), \ \ \ cos(2\pi z), \ \ \ \wp(z,1,\tau), \ \ \  \wp^{\prime }(z,1,\omega),\ \ \ \wp^2(z,1,i), 
\ \ \ \wp^{\prime 2}(z,1,\omega),$$
where $\wp=\wp(z,\omega_1,\omega_2)$ denotes	the
Weierstrass functions with periods $\omega_1,$ $\omega_2$ (see \cite{dh}, \cite{mil2}).

In terms of orbifolds, the Ritt theorem can be formulated as follows (\cite{e2}). 

\bt Let $A$ and $C$ be commuting rational functions of degree at least two having no iterate in common. 
Then there exists an orbifold  $\f O=(\C\P^1,\nu)$ with $\chi(\f O)=0$ such that 
$A:\f O\rightarrow \f O$ and 
$C:\f O \rightarrow \f O$ are covering maps between orbifolds. \qed
\et

If $\f O=(\C\P^1,\nu)$ is an orbifold with $\chi(\f O)=0$, and $f$ is a  
rational function such that $f:\f O \rightarrow \f O$  is a covering map  between orbifolds, then $\t{R}=\C$, and $f$ 
lifts to an affine map $\phi=a z+b,$  $a,b\in \C,$ which makes
the diagram 
\be \l{ida}
\begin{CD}
\C @>\phi=a z+b >> \C\\
@VV\theta_{\f O}V @VV\theta_{\f O}V\\ 
\f O @>f >> \f O, 
\end{CD}
\ee
commutative.
Thus,  on one hand, the Ritt theorem reduces describing pairs  of commuting rational functions $A$ and $C$ having no iterate 
in common to describing pairs of affine maps 
$\phi$ and $\psi$ that map any orbit of some group $\Gamma$ from  list \eqref{lis} to another orbit 
and satisfy the equality 
$$\phi\circ \psi= g\circ \psi\circ \phi$$
for some $g\in \Gamma$. On the other hand,  the Ritt theorem imposes restrictions on possible 
ramifications of $A$ and $C$ resulting from the definition of covering map \eqref{us} and list \eqref{list}.
Note that if  $\f O=(\C\P^1,\nu)$ is an orbifold and $f:\f O \rightarrow \f O$ is a covering map of degree at least two, then \eqref{rhor} implies that $\chi(\f O)=0$. In particular, the condition $\chi(\f O)=0$ in the formulation of the Ritt theorem is actually redundant.

If $\nu(\f O)=(\infty,\infty)$, then   
any rational function $f$ of degree at least two such that $f:\f O \rightarrow \f O$ is a  covering map between orbifolds is conjugate to $z^{\pm n}$, $n\geq 2,$ while 
if $\nu(\f O)=(2,2,\infty)$, then any such  function is conjugate to $\pm T_n,$ $n\geq 2.$ 
Rational functions $f$ of degree at least two such that $f:\f O \rightarrow \f O$ is a covering map for an orbifold $\f O$  whose signature is $(3,3,3),$ $(2,4,4),$ $(2,3,6),$ or $(2,2,2,2)$ are 
 called {\it Latt\`es maps}. Such rational functions possess a number of remarkable  features  (see \cite{mil2}, \cite{gl}).

\vskip 0.2cm

In this paper  all considered orbifolds (except for universal coverings) will be defined on $\C\P^1,$ and  
we simply will write $\f O$ instead of $\f O=(\C\P^1,\nu)$.
The following statement describes compositional properties of rational functions $C$ that are self-covering maps  $C:\f O\rightarrow \f O$ (cf. \cite{gl}, Theorem 4.1).

\bl \l{meba} Let $\f O$ be an orbifold and $C$ a rational function such that \linebreak $C:\f O\rightarrow \f O$ is a covering map between orbifolds. Assume that $C=X\circ Y$ is a 
decomposition of $C$ into a composition of rational functions.
Then there exists an orbifold 
$\f O^*$ with $\nu(\f O^*)=\nu(\f O)$ such that
$ Y:\f O \rightarrow \f O^*$ and  
$X:\f O^* \rightarrow \f O$
are covering maps between orbifolds.
\el 
\pr 
Since 
\be \l{leg} \nu((X\circ Y)(z))=\nu(z)\deg_z(X\circ Y)=\nu(z)\deg_zY\deg_{Y(z)}X\ee 
and the value $(X\circ Y)(z)$ depends only on the value $Y(z)$, 
defining for $z\in \C\P^1$ the value $\nu^*(z)$ by the formula  
$$\nu^*(z)=\nu(z')\deg_{z'}Y, $$ where $z'$ is any point  such that $Y(z')=z,$ 
we obtain a well-defined orbifold $\f O^*$ such that $ Y:\f O \rightarrow \f O^*$ and  
$X:\f O^* \rightarrow \f O$ are covering maps. Moreover, 
applying formula \eqref{rhor} to any of these maps we see that $\chi(\f O^*)=0$. Finally, it is not hard  to prove 
that $\nu(\f O^*)=\nu(\f O)$. Indeed, if $\nu(\f O)=(\infty,\infty)$, then  \eqref{leg} implies easily that $\f O^*$ has exactly two points with ramification $\infty$. Therefore,  since $\f O^*$ belongs to list \eqref{list}, the equality $\nu(\f O^*)=(\infty,\infty)$ holds. 
 Similarly, we obtain that  if $\nu(\f O)=(2,2,\infty)$, then $\nu(\f O^*)=(2,2,\infty)$. Assume now that 
$\nu(\f O)=(2,3,6)$. Since $X:\f O^* \rightarrow \f O$ is a covering map, it follows from \eqref{us} that 
$$\nu^*(z)\mid \nu(X(z)),\ \ \  z\in \C\P^1,$$  implying  that either $\nu(\f O^*)=(2,3,6)$, or 
$\nu(\f O^*)=(3,3,3)$, or $\nu(\f O^*)=(2,2,2,2)$. However, in the last two cases 
$Y:\f O \rightarrow \f O^*$ cannot be a covering map,
since $$\nu(z)\mid \nu^*(Y(z)),\ \ \  z\in \C\P^1.$$
The rest of the cases are considered similarly. 
\qed

\vskip 0.2cm

Let us list  several  properties of Latt\`es maps used in the following. 
First, if $f$ is a Latt\`es map, then an orbifold $\f O$ 
such that $f:\f O\rightarrow \f O$ is a covering map, is defined in a unique way by dynamical properties of $f$ 
(see \cite{mil2} and also \cite{gl}, Theorem 6.1). 
We will use the notation $\f O_f$ for this orbifold and the notation   $l=l(f)$ 
for  the least common multiple of numbers in the signature of $\f O_f$.
Secondly, although the functions $\theta_{\f O}$ and $\phi$ in diagram \eqref{ida} are not defined in a unique way by $f$, the number $a^l$ depends on $f$ only, and the numbers $a$ and $\deg f$ 
are related by the equality \be \l{ryba} \deg f=\vert a\vert^2\ee 
(see \cite{mil2}, Lemma 5.1).
Thirdly, if $f$ satisfies \eqref{ida} and $z\in \C\P^1$ is a fixed point of $f$, then  
the multiplier of $f$ at  $z$ is given by the formula  
\be \l{for} \mu=(\omega a)^{\nu(z)},\ee
where $\omega$ is some $l$th root of unity, and $\nu$ is the ramification function for $\f O_f$ (see \cite{mil2}, Corollary 3.9). 

\vskip 0.2cm

Finally, we need the following rigidity  property of Latt\`es  maps 
which states, roughly speaking, that if $l\geq 3$, then for fixed  $a^l$ there exist at most two conjugacy classes of rational functions $f$ 
which make diagram \eqref{ida} commutative, and these classes can be distinguished by their dynamical properties 
(see \cite{mil2}, Theorem 5.2). 

\bt \l{xork}  Let  $f$ be a Latt\`es map with $l=l(f)\geq 3.$
Then the conjugacy class of $f$
is completely determined by the numbers $l$ and $a^l$ together with the information as to whether $f$ does or does not have a fixed point of multiplier $\mu=a^l.$ \qed
\et 
Note that  in view of formula \eqref{for} the property of $f$ to have a fixed point of multiplier $\mu=a^l$ is equivalent to the following property:

\vskip 0.2cm

($\star$) {\it there exists a fixed point $z$ of $f$ with $\nu(z)=l$.} 
\vskip 0.2cm 

Theorem \ref{xork} results in the following statement. 

\bt \l{medvi} Let $A$ and $B$ be mutually semiconjugate rational functions of degree at least two, and $X$, $Y$  rational functions such that diagram \eqref{en} commutes. Assume that there exists an orbifold $\f O$ with signature distinct from $(2,2,2,2)$  
such that $A:\f O\rightarrow \f O$ and $X\circ Y:\f O \rightarrow \f O$ are covering maps between orbifolds. Then $B$ is conjugate to $A$. 
\et
\pr By Lemma \ref{meba}, there exists 
an orbifold $\f O^*$  with  $\nu(\f O^*)=\nu(\f O)$ such that 
 $ Y:\f O \rightarrow \f O^*$ and  
$X:\f O^* \rightarrow \f O$
are covering maps between orbifolds. Furthermore, 
 since  $$\nu(\f O^*)=\nu(\f O)\neq (2,2,2,2),$$  
changing in diagram \eqref{en} the function $Y$ to the function $\mu\circ Y$, the function $X$ to the function  $X\circ \mu^{-1}$, 
 and the function $B$ to the function $\mu^{-1}\circ B\circ \mu$ 
for convenient 
$\mu\in Aut(\C\P^1)$, 
we may assume 
that $\f O^*=\f O.$ 

If $\nu(\f O)=(\infty,\infty)$, then
without loss 
of generality we may assume that $$\nu(0)=\infty, \ \ \ \nu(\infty)=\infty,$$ implying that 
$$A=az^n, \ \ \ Y=bz^m,$$ where $a,b\in \C$ and  $n,m\in \Z.$ 
It follows now from 
the equality $B\circ Y=Y\circ A$ that $B=a^mb^{1-n}z^n$. Thus, in this case $A$ and $B$ are conjugate.

Similarly, if $\nu(\f O)=(2,2,\infty)$ and 
$$\nu(1)=2, \ \ \ \nu(-1)=2, \ \ \ \nu(\infty)=\infty,$$
then $$A=\pm T_n, \ \ \  Y=\pm T_{m_1},$$ implying that 
$B=\pm T_n$.
 However, in this case a further investigation is needed, 
since the functions $T_n$ and $-T_n$ are conjugate for even $n$, but  not conjugate 
for odd. 
To finish the proof, we observe that 
the equality  
\be \l{wp} -T_n\circ \pm T_m=\pm T_m\circ T_n\ee for odd $n$ is impossible. Thus,  
 if $A=T_n$, then $B=T_n$. 
In turn, this implies that  
if $A=-T_n$, then $B=-T_n$,  
for otherwise the lower square in \eqref{en} would provide a 
solution of \eqref{wp}.

Finally, assume that $\nu(\f O)$ is $(2,4,4),$ $(3,3,3)$, or $(2,3,6)$. Let us complete diagram \eqref{en} to the diagram  
\be \l{kro}
\begin{CD}
\C @>\phi=a z+b >> \C\\
@VV\theta_{\f O}V @VV\theta_{\f O}V\\ 
\C\P^1 @>A>> \C\P^1\\
@VV Y V @VV Y V\\ 
\C\P^1 @>B >> \C\P^1\\ 
@VV X V @VV X V\\ 
\C\P^1 @>A >> \C\P^1\ . 
\end{CD}
\ee
Since $\theta_{\f O}: \C\rightarrow \f O$  and $Y:\f O \rightarrow \f O$ are covering maps, their composition \linebreak
$Y\circ  \theta_{\f O}: \C\rightarrow \f O$ is also a covering map
(here $\C$ stands for the orbifold $(\C, \nu)$ with $\nu \equiv 1$).
Thus,  $Y\circ  \theta_{\f O}$ along with $\theta_{\f O}$ is a universal covering of $\f O,$ 
implying by the chain rule that 
$B:\f O \rightarrow \f O$ is a covering map. 
Moreover, \eqref{kro} implies that $A$ and $B$  
have the same invariant  $a^l$. 
Therefore, by Theorem \ref{xork}, 
it is enough to show that  property ($\star$) holds 
for $A$ if and only if it holds for $B$.

Consider the semiconjugacy in  the upper square in \eqref{en}. Clearly, $Y$ maps fixed points of $A$  to fixed points of 
$B$. Therefore, since  the equality
$$\nu(Y(z))=\nu(z)\deg_zY$$ implies that $\nu(Y(z))=l$ whenever $\nu(z)=l$,
 if  property ($\star$) holds for $A$, then it holds for $B$. Moreover, arguing as 
 in the case $\nu(\f O)=(2,2,\infty)$, we conclude  
that if the property ($\star$) does not hold 
for $A$, then it does not hold for $B$. \qed

\end{section}

\begin{section}{Equivalence  and semiconjugacy}
Let $A$ be a rational function. Recall that for any decomposition $A=U\circ V$ of $A$ into a composition 
of rational functions the 
rational function $\t A=V\circ U$ is called an {\it elementary transformation} of $A$, and rational functions $A$ and $B$ are called {\it equivalent} if there exists 
a chain of elementary transformations between $A$ and $B$. 
Since for any M\"obius transformation  $\mu$ the equality
$$A=(A\circ \mu)\circ \mu^{-1}$$ holds, each equivalence class $[A]$ is a union of conjugacy classes. 
Thus, like the mutual semiconjugacy relation, the relation $\sim$ is a 
weaker form of the classical conjugacy relation. Moreover, equivalent rational functions have similar dynamic characteristics. 
To make the last statement precise, recall that the {\it multiplier spectrum} of a rational function $A$ of degree $d$ is a function which assigns to each $s\geq 1$ the unordered list of multipliers at all $d^s+1$ fixed points of $A^{\circ s}$ taken with appropriate multiplicity. 
Two rational functions  are called {\it isospectral} if they have the same   multiplier spectrum.
In this notation, the following statement is true  (see \cite{rec}, Corollary 2.1).

\bl \l{uv+} Let $A$ and $B$ be rational functions such that $A\sim B$. Then $A$ and $B$ are isospectral. \qed
\el

Lemma \ref{uv+}  has two implications. On one hand, it permits to 
conclude that two functions are {\it not} equivalent if they have different   multiplier spectrum. On the other hand, by the fundamental result of McMullen (\cite{Mc}),  
the conjugacy class of any rational function $A$ that is not {\it a flexible Latt\`es map} (see e.g. \cite{mil2} or \cite{sildyn} for the definition) is defined up to finitely many choices by its
multiplier spectrum. Thus, Lemma \ref{uv+} implies that for any function $A$ that is not a flexible Latt\`es map
the number of  conjugacy classes in the equivalence class $[A]$ is finite. More precisely, the following statement 
holds (see \cite{rec}, Theorem 1.1).

\bt \l{t1} Let $A$  be a rational function. Then its  equivalence class $[A]$ 
contains infinitely many conjugacy classes if and only if 
$A$ is a flexible Latt\`es map. \qed
\et

Note that there exists no  {\it absolute} bound for the number of conjugacy classes in $[A]$, and one can construct  rational functions $A$ of degree $n$ for which $[A]$ contains  $\approx\log_2n$ conjugacy classes (see \cite{semi}, p. 1241). On the other hand, although the proof of the McMullen theorem is non-effective, Theorem \ref{t1} can be deduced from effective results of the paper \cite{pnew}, implying that the number of conjugacy classes in $[A]$ can be bounded in terms of degree of $A$ only.  

Finally, we mention that to our best knowledge only three types of examples 
of isospectral rational functions are known: flexible or not flexible Latt\`es maps (see \cite{Mc}, \cite{mil2}, \cite{sildyn}), and 
equivalent functions. So, the following question is of great interest.

\bp Do there exist isospectral rational functions that are neither Latt\`es maps nor equivalent ?
\ep

It was already mentioned in the introduction that the equivalence  $\sim$ is closely related to 
the semiconjugacy. Moreover, using elementary transformations one can reduce any solution of \eqref{i1}
to a so-called primitive solution.  
We say that a solution $A,X,B$  
of functional equation \eqref {i1} is   {\it primitive} if $$ \C(X,B)=\C(z),$$
	that is if the functions $X$ and $B$ generate the whole field of rational functions. 
It was shown in \cite{semi} (see also \cite{arn}), that for any primitive solution of \eqref{i1}
there exist orbifolds $\f O_1$ and $\f O_2$ 
such that $A:\f O_1\rightarrow \f O_1$, $B:\f O_2\rightarrow \f O_2$, and $X:\f O_1\rightarrow \f O_2$ are {\it minimal holomorphic maps} between orbifolds (see  \cite{semi} for the definition).
This condition generalizes the condition provided by the Ritt theorem, and implies strong restrictions on a possible  form of $A$, $B$ and $X$. 

In what follows, we will not use the description of primitive solutions given in \cite{semi}. However, we will need
the following reduction: for an arbitrary solution  $A,X,B$   of \eqref{i1} there exists a decomposition $X=X_0\circ W$ and a rational function $B_0\sim B$
such that the diagram 
\be \l{kor}
\begin{CD} 
\C\P^1 @> B>>\C\P^1 \\ 
@V W  VV @VV W   V\\ 
 \C\P^1 @> B_0>> \C\P^1
 \\ 
@V {X_0}  VV @VV {X_0}   V\\ 
 \C\P^1 @> A>> \C\P^1, 
\end{CD} 
\ee
commutes and
$A,X_0,B_0$ is a primitive solution of \eqref{i1}.  To see that this is true, we observe that 
if $A,X,B$ is a primitive solution, then we can set $W=z,$ $X_0=X,$ $B_0=B.$ On the other hand, if the solution  $A,X,B$ is not primitive, then by the L\"uroth theorem
 there exists a rational function $W$ of degree greater than one such that
 $\C(X,B)= \C(W)$ and the equalities
$$ X=X'\circ W,\ \ \  B=B'\circ W $$ hold for some rational functions $X'$ and $B'$ with 
$\C(X',B')=\C(z).$
Clearly, the diagram 
$$
\begin{CD} 
\C\P^1 @> B>>\C\P^1 \\ 
@V W  VV @VV W   V\\ 
 \C\P^1 @> W\circ B'>> \C\P^1
 \\ 
@V {X'}  VV @VV {X'}   V\\ 
 \C\P^1 @> A>> \C\P^1, 
\end{CD} 
$$
commutes. Thus, if the solution  $A, X',W\circ B'$  of \e{i1} is  primitive, we are done.
Otherwise, we can apply the above transformation to this solution. Since
$\deg X'<\deg X$, it is clear that after a finite number of steps we will obtain
required functions $X_0,B_0, W$  (see \cite{gl}, Sec\-tion 3, for more details).

In addition to the above reduction, to prove Theorem \ref{main} we need the following two lemmas (see \cite{rev}, Lemma 2.4 and Lemma 2.5).

\bl \l{asd} A solution $A,X,B$ of \eqref{i1} is primitive if and only if the algebraic 
curve 
\be \l{cur} A(x)-X(y)=0\ee is irreducible. \qed
\el

\bl \l{primsol}
Let $A,X,B$ be a primitive solution of \eqref{i1}. Then for any $s\geq 1$ the triple $A^{\circ s},X,B^{\circ s}$ is also a primitive solution of \eqref{i1}. \qed
\el

\noindent{\it Proof of Theorem \ref{main}}. Let $A$, $B$ be
mutually semiconjugate rational functions, and $X$, $Y$ corresponding rational functions which make diagram \eqref{en} commutative. Then by the Ritt theorem, either there exists an orbifold  $\f O$ with $\chi(\f O)=0$ such that 
$A:\f O\rightarrow \f O$ and 
$X\circ Y:\f O \rightarrow \f O$
are covering maps between orbifolds, or
there exist  $s,k\geq 1$ such that \be \l{medv} A^{\circ s}=(X\circ Y)^{\circ k}.\ee
In the first case, the statement of the theorem follows from  Theorem \ref{medvi}. On the other hand, 
to prove the theorem in the second case it is enough to show that 
in diagram \eqref{kor}, constructed for $A,X,B$ 
from the lower square in \eqref{en}, the equality   
$\deg X_0=1$ holds   (cf. \cite{rev}, Theorem 2.5). Indeed, in this case  $B_0$ is conjugate to $A$, and hence 
$$B\sim B_0\sim A.$$ 

Assume to the contrary that $\deg X_0\geq 2.$ Set $$F=W\circ Y\circ  (X\circ Y)^{\circ k-1},$$
where $k$ is defined by \eqref{medv}. 
Then
$$A^{\circ s}=X_0\circ F$$ by \eqref{medv},  implying that 
the curve 
$$F(x)-y=0$$ is a component of the curve
\be \l{cucu} A^{\circ s}(x)-X_0(y)=0.\ee Moreover, since $\deg X_0>1,$ this component is proper. Therefore, the triple $A^{\circ s},X_0,B_0^{\circ s}$ is not a primitive solution of \eqref{i1} by Lemma \ref{cur}. On the other hand, this triple must be 
 a primitive solution  by Lemma \ref{primsol}. The contradiction obtained shows that 
$\deg X_0=1$. 
 \qed

\end{section}

\begin{section}{Case of signature $(2,2,2,2)$} 
We recall that any  Latt\`es map with invariant $l$ equal to 2 is conjugate to a rational function  
$f$ such that the diagram 
\be \l{avo}
\begin{CD}
\C @>\phi=az+b >> \C\\
@VV \wp_L V @VV \wp_L V\\ 
\f O @>f >> \f O 
\end{CD}
\ee commutes for some   lattice  $L$ of rank two in $\C$ and affine map $\phi$. The Weierstrass function $\wp_L$ is the universal covering of $\f O_f,$
and the corresponding group $\Gamma_{\f O}$ is generated by translations by elements of  $L$  and the transformation  $z\rightarrow -z.$ Furthermore, the 
function $\phi=az+b$ in \eqref{avo} maps
any orbit of $\Gamma_{\f O}$ to another orbit, 
implying that $aL\subset L$ (see e.g. \cite{mil2}, Lemma 5.1). For most lattices $L$ the condition $aL\subset L$
implies that $a\in \Z$. In particular,  for such $L$ the degree of $f$ in \eqref{avo} is a perfect square by formula \eqref{ryba}.  
Lattices for which there exists a non-integer $a$ satisfying $aL\subset L$ are called lattices with {\it complex multiplication}. 

For an integer $n\geq 2$ and a lattice $L$ we define a Latt\`es map $A_{n,L}$ by the 
commutative diagram 
\be \l{lev} 
\begin{CD}
\C @>\phi=n z >> \C\\
@VV \wp_L V @VV \wp_L V\\ 
\f O @>A_{n,L} >> \f O. 
\end{CD}
\ee
Clearly, the rational functions $A_{n,L}$ and $A_{m,L}$ commute for any $n,m\geq 2$. 
Let $L'$ be a lattice satisfying \be \l{sat} L\subset L'\subset L/n.\ee 
For example, if $L=\langle\omega_1,\omega_2\rangle$, we can set $L'=\langle \omega_1,\frac{\omega_2}{n}\rangle$.
With such $L'$ we can associate  
a functional decomposition \be \l{hol} A_{n,L}=X_{L'}\circ Y_{L'}\ee as follows.
Since any even doubly periodic meromorphic function with period lattice $L$ is a rational function in $\wp_L,$ it follows from 
\eqref{sat} that there exist rational functions 
$X$, $Y$, $F$ such that 
$$ \wp_{L/n}=X \circ \wp_{L'}, \ \ \ \  \wp_{L'}=Y \circ \wp_{L},\ \ \ \ \wp_{L/n}=F \circ \wp_{L},$$ and  
it is clear that $F=X\circ Y.$
On the other hand, since 
$$ \wp_{L/n}(z/n)=n^2\wp_L(z),$$ we have:
$$\wp_{L/n}=n^2\wp_L(nz)=n^2A_{n,L}\circ \wp_L.$$ Therefore, $F=n^2A_{n,L}$, and hence \eqref{hol} holds
for $X_{L'}=X/n^2$ and $Y_{L'}=Y.$
Another way to obtain decomposition \eqref{hol} is to consider the projections of the isogeny 
$\C/L\rightarrow \C/L'$ and its dual  (see \cite{rec}, Section 3).
Note that  
\be \l{eg} \deg X_{L'}=[L/n:L'], \ \ \ \ \deg Y_{L'}=[L':L],\ee
and both these numbers are greater than one since  $L'$ is distinct from $L$ and $L/n.$ Finally, it is clear that
\be \l{eg2} \deg X_{L'}\cdot\deg Y_{L'}=n^2.\ee

We now observe that since the diagram 
$$
\begin{CD}
\C\P^1 @>mz>> \C\P^1\\
@VV \wp_L V @VV \wp_L V\\ 
\C\P^1 @>A_{m,L} >> \C\P^1\\ 
@VV A_{n,L} V @VV A_{n,L} V\\ 
\C\P^1 @>A_{m,L} >> \C\P^1\ 
\end{CD}
$$
commutes, it follows from the equalities  \eqref{hol} and $\wp_{L'}=Y_{L'} \circ \wp_{L}$  
that the diagram 
\be \l{oxb}
\begin{CD}
\C\P^1 @>A_{m,L}>> \C\P^1\\
@VV Y_{L'} V @VV Y_{L'} V\\ 
\C\P^1 @> A_{m,L'}>> \C\P^1\\ 
@VV X_{L'} V @VV X_{L'} V\\ 
\C\P^1 @>A_{m,L} >> \C\P^1\ 
\end{CD}
\ee
also commutes. Thus, $A_{m,L'}$ and $A_{m,L}$ are mutually semiconjugate.
Since for any lattice $L$ there exist lattices $L'$ satisfying \eqref{sat}, we obtain in this way a large class of examples of mutually semiconjugate rational functions, and we will show below that at least some of these functions  are not equivalent. 

Note that 
the functions 
$A_{m,L'}$ and $A_{m,L}$ are  isospectral for any lattice $L$, since the multiplier spectrum of $A_{m,L}$ depends only on $m$ (see e.g. \cite{sildyn}, Proposition 6.52(b)). Thus, we cannot use Lemma \ref{uv+} to prove that $A_{m,L'}\not\sim A_{m,L}$. Instead, we  use
the following observation.

\bl \l{esd} Let $A$ be a Latt\`es map of degree $d$ with $l(A)=2$, and $B$  a rational function such that $B\sim A$. Then $B$ is a Latt\`es map with $l(B)=2$, and  there exists  a rational function $T$ whose degree divides $d^k$, $k\geq 0,$ such that  $T: \f O_A\rightarrow \f O_B$  is a covering map.
\el
\pr 
If $B$ is an elementary transformation of $A$, that is 
$A=U\circ V$ 
and  $B=V\circ U$ for some rational functions  $U$ and $V$,  
then by Lemma \ref{meba}  there exists an orbifold $\f O'$ with $\nu(\f O')=(2,2,2,2)$ such that 
$$V:\f O \rightarrow \f O', \ \ \ U:\f O' \rightarrow \f O$$
 are covering maps.
Therefore,  $B: \f O'\rightarrow \f O'$ is a covering map, and hence $B$ is a Latt\`es map with $\f O_B=\f O'$ and $l(B)=2$.    
Moreover, the map $T=V$ satisfies the  requirements of the lemma since $\deg V$ is a divisor of $d$. 
Since any $B\sim A$ is obtained from $A$ by a chain of elementary transformations, and elementary transformations do not change the degree,
using the above reasoning recursively and composing corresponding functions $V$, 
we obtain a rational function $T$ with the required properties. \qed

\bt \l{lux} Let $L$ be a lattice without complex multiplication, and $n$, $m$ distinct primes.
Then for any lattice  $L'$  satisfying  $ L\subset L'\subset L/n$ the functions 
$A_{m,L'}$ and $A_{m,L}$ are mutually semiconjugate but non-equivalent.
\et
\pr  
Assume that ${A_{m,L}}\sim A_{m,L'}$, and let $T: \f O_{A_{m,L}}\rightarrow \f O_{{A_{m,L'}}}$be  a covering map 
between orbifolds provided by Lemma \ref{esd}. Then $\deg T=m^k$, $k\geq 0,$ since $\deg A_{m,L}=m^2$ and $m$ is a prime.  
Applying Lemma \ref{meba} to decomposition \eqref{hol}, we conclude that there exists an orbifold $\f O^*$ with $\nu(\f O^*)=(2,2,2,2)$ 
such that 
$$Y_{L'}:\f O \rightarrow \f O^*, \ \ \ X_{L'}:\f O^* \rightarrow \f O$$ 
are covering maps, and as in the proof of Theorem \ref{medvi} we see that the map 
$\wp_{L'}=Y_{L'} \circ \wp_{L}$ is the universal covering of $\f O^*$. Since $\wp_{L'}$ is also 
the universal covering of $\f O_{A_{m,L'}}$, this implies that 
 \be \l{az} \f O^*=\f O_{A_{m,L'}}.\ee  Thus,
$X_L':\f O_{A_{m,L'}} \rightarrow \f O_{A_{m,L}}$ is  a covering map, and hence  the composition
$$X_L'\circ T:\, \f O_{A_{m,L}}\rightarrow \f O_{A_{m,L}}$$ is also a covering map.
Since by assumption $L$ is a lattice without complex multiplication, the number $\deg (X_L'\circ T)$ must be  a perfect square. 
On the other hand, since $n$ is a prime, it follows from \eqref{sat} and \eqref{eg2} that 
$$\deg X_L'=[L/n:L']=n,$$
implying that $$\deg (X_L'\circ T)=nm^k,\ \ \ k\geq 0.$$ Thus, since $n>1$ and $\gcd(n,m)=1$, the number $nm^k$ cannot be a 
perfect square.
The contradiction obtained finishes the proof. \qed

\end{section} 

\vskip 0.2cm

\noindent{\it Acknowledgments.}
This research was supported by the ISF Grant No. 1432/18.

\bibliographystyle{amsplain}

\end{document}